\documentclass[11pt]{article}
\usepackage{graphicx}
\usepackage{color}
\usepackage{amssymb,latexsym,amsmath}
\usepackage{epstopdf}
\newtheorem{lemma}{Lemma}[section]
\newtheorem{theo}[lemma]{Theorem}

\newtheorem{coro}[lemma]{Corollary}
\newtheorem{remark}[lemma]{Remark}

\newcommand{\proof}{\noindent{\em Proof: }}

\newcommand{\forme}[1]{}

\def\wbull{\hfill\vrule height .9ex width .8ex depth -.1ex}

\begin{document}

\date{\today}
\title{The Terwilliger polynomial of a $Q$-polynomial\\ distance-regular graph and its application\\
to the pseudo-partition graphs}
\author{{\bf Alexander L. Gavrilyuk}\\
Research Center for Pure and Applied Mathematics,
Graduate School of Information Sciences,\\
Tohoku University, Sendai 980-8579, Japan\\
and\\
N.N. Krasovsky Institute of Mathematics and Mechanics UB RAS,\\ 
Kovalevskaya str., 16, Ekaterinburg 620990, Russia\\
e-mail: alexander.gavriliouk@gmail.com\\
\\
{\bf Jack H. Koolen}\\School of Mathematical Sciences\\ University of Science and Technology of China, 
Hefei, 230026, Anhui, PR China\\
e-mail: koolen@ustc.edu.cn}
\maketitle

\begin{abstract}
Let $\Gamma$ be a $Q$-polynomial distance-regular graph with diameter at least $3$. 
Terwilliger (1993) implicitly showed that there exists a polynomial, say $T(\lambda)\in \mathbb{C}[\lambda]$, 
of degree $4$ depending only on the intersection numbers of $\Gamma$ and such that $T(\eta)\geq 0$ holds
for any non-principal eigenvalue $\eta$ of the local graph $\Gamma(x)$ for any vertex $x\in V(\Gamma)$. 

We call $T(\lambda)$ the Terwilliger polynomial of $\Gamma$. In this paper, we give an explicit 
formula for $T(\lambda)$ in terms of the intersection numbers of $\Gamma$ and its dual eigenvalues. 
We then apply this polynomial to show that all pseudo-partition graphs with diameter at least $3$ are known.
\end{abstract}

\section{Introduction}

Let $\Gamma$ denote a distance-regular graph with diameter $D\geq 3$. 
(For definitions and notations see Section \ref{SectDefinitions}.)
Pick any 3-tuple $xyz$ of vertices of $\Gamma$ such that $y$ and $z$ are 
neighbours of $x$. 
Let $[i,j,k]$ denote the number of vertices $u$ of $\Gamma$ 
such that $u$ is at distances $i$ from $x$, $j$ from $y$, and $k$ from $z$, 
respectively. In general, $[i,j,k]$ depends on the choice of $x,y$ and $z$ as well as $i,j,k$. 

Assuming that $\Gamma$ has a $Q$-polynomial structure, Terwilliger showed in \cite{KiteFree} that 
if $x,y,z$ are mutually adjacent then
\begin{equation*}
[i,i-1,i-1]=\kappa_{i}[2,1,1]+\tau_{i}
\end{equation*}
holds, where $\kappa_{i},\tau_{i}$ are real 
scalars that do not depend on $x,y,z$. 

Inspired by his proof, we show in Section \ref{spearSection} under the same assumption on $\Gamma$ 
that 
\begin{equation}\label{spear1}
[i,i+1,i+1]=\sigma_{i}[1,2,2]+\rho_{i,\delta},
\end{equation}
where $\delta$ is the distance between $y$ and $z$, $y\ne z$, i.e. $\delta\in \{1,2\}$, 
and $\sigma_{i},\rho_{i,\delta}$ are real 
scalars that do not depend on $x,y$, and $z$. In fact, 
$\sigma_{i}$ and $\rho_{i,\delta}$ are rational expressions 
in the dual eigenvalues and intersection numbers of $\Gamma$.

Terwilliger (see 'Lecture note on Terwilliger algebra' edited by Suzuki, \cite{SN}) showed 
that, for $i=2,3,\ldots,D-1$, there exists 
a polynomial $T_i(\lambda)\in \mathbb{C}[\lambda]$ of degree $4$ such that for any $i$, 
any vertex $x\in V(\Gamma)$, and any non-principal eigenvalue $\eta$ 
of the local graph $\Gamma(x)$ one has $T_i(\eta)\geq 0$.

In \cite{SN} it was claimed that $T_i(\lambda)$ does not depend on $i$ up to a scalar multiple. 
We have no proof for this, but see Remarks \ref{Tclaimremark} and \ref{Tclaimdiscuss}. 
We call $T_i(\lambda)$ the Terwilliger polynomial.

In Section 4 we will give a proof of the existence of $T_i(\lambda)$ and calculate the polynomial 
$T_2(\lambda)$ explicitly for any $Q$-polynomial distance-regiular graph. 

Much attention has been paid to a project of classification 
of all $Q$-polynomial distance-regular graphs with large diameter.
Bannai and Ito's interpretation of the Leonard theorem says that the intersection numbers 
of a $Q$-polynomial distance-regular graph has one of seven types: 
$1$, $1A$, $2$, $2A$, $2B$, $2C$, and $3$, see \cite[Theorem 5.1]{BI} or \cite{SubAlgPaper}.

In the survey paper \cite{SurveyDRG} by Van Dam, Koolen and Tanaka the current 
status of the classification of the $Q$-polynomial distance-regular graphs is reported.

In Section 5 we will focus on the $Q$-polynomial distance-regular graphs of type $2$.
Using the Terwilliger polynomial and under some additional restrictions on parameters 
(which cover all the known examples), 
we classify the $Q$-polynomial distance-regular graphs of type $2$, including some of the pseudo-partition graphs. 
Thus we settle the uniqueness problem for all open cases of pseudo-partition graphs 
with diameter at least $3$.

In a forthcoming paper, we will apply the Terwilliger polynomial to show 
the uniqueness of the Grassmann graphs $J_2(2D,D)$ with odd diameter $D$ 
and the bilinear forms graphs $H_2(D,D)$ for all $D\geq 3$.

\section{Definitions and preliminaries}\label{SectDefinitions}

All the graphs considered in this paper are finite, undirected and simple. 
Suppose that $\Gamma$ is a connected graph with vertex set $V(\Gamma)$ and 
edge set $E(\Gamma)$, where $E(\Gamma)$ consists of unordered pairs of  
adjacent vertices. The distance $d(x,y)$ between any two vertices $x,y$ of
$\Gamma$ is the length of a shortest path connecting $x$ and $y$ in $\Gamma$. 

For a subset $X$ of the vertex set of $\Gamma$, we will also 
write $X$ for the subgraph of $\Gamma$ induced by $X$.
For a vertex $x\in V(\Gamma)$, define $\Gamma_i(x)$ 
to be the set of vertices which are at distance precisely $i$ from $x$ ($0\leq i\leq D$),
where $D:={\rm max}\{d(x,y)\mid x,y\in V(\Gamma)\}$ is the {\it diameter} of $\Gamma$. 
In addition, define $\Gamma_{-1}(x)=\Gamma_{D+1}(x)=\emptyset$. 
The subgraph induced by $\Gamma_1(x)$ is called the {\it neighborhood} or the {\it local graph} of a vertex $x$.
We write $\Gamma(x)$ instead of $\Gamma_1(x)$ for short, and we denote $x\sim_{\Gamma} y$ or simply $x\sim y$ 
if two vertices $x$ and $y$ are adjacent in $\Gamma$. 
For a graph $G$, a graph $\Gamma$ is called {\it locally} $G$ if any local graph of $\Gamma$ is isomorphic to $G$.

The {\it eigenvalues} of a graph are the eigenvalues of its adjacency matrix.
If, for some eigenvalue $\eta$ of $\Gamma$, its eigenspace contains 
a vector orthogonal to the all ones vector, we say the eigenvalue $\eta$ is {\it non-principal}.
If $\Gamma$ is regular with valency $k$ then all its eigenvalues are non-principal 
unless the graph is connected and then the only eigenvalue that is principal 
is its valency $k$.

A connected graph $\Gamma$ with diameter $D$ is called {\it distance-regular} 
if there exist integers $b_{i-1}$, $c_{i}$ $(1\leq i\leq D$) such that, for any two 
vertices $x,y\in V(\Gamma)$ with $d(x,y)=i$, there are precisely $c_i$ neighbors 
of $y$ in $\Gamma_{i-1}(x)$ and $b_i$ neighbors of $y$ in $\Gamma_{i+1}(x)$. 
In particular, any distance-regular graph is regular with valency $k:=b_0$. 
We define $a_i:=k-b_i-c_i$ for notational convenience and note that 
$a_i=|\Gamma(y)\cap \Gamma_i(x)|$ holds for any two vertices $x,y$ with 
$d(x,y)=i$ ($1\leq i\leq D$). 
The array $\{b_0,b_1,\ldots,b_{D-1};c_1,c_2,\ldots,c_D\}$ is 
called the {\it intersection array} of the distance-regular graph $\Gamma$.

A distance-regular graph with diameter 2 is called a {\it strongly regular} graph. 
We say that a strongly regular graph $\Gamma$ has parameters $(v,k,\lambda,\mu)$, 
if $v=|V(\Gamma)|$, $k$ is its valency, $\lambda:=a_1$, and $\mu:=c_2$.

If a graph $\Gamma$ is distance-regular then, for all integers $h,i,j$ 
($0\leq h,i,j\leq D$), and all vertices $x,y\in V(\Gamma)$ with $d(x,y)=h$, 
the number 
$$p^h_{ij}:=|\{z\in V(\Gamma)\mid d(x,z)=i,~d(y,z)=j\}|$$
does not depend on the choice of $x,y$. The numbers $p^h_{ij}$ 
are called the {\it intersection numbers} of $\Gamma$.
Note that $c_i=p^i_{1i-1}$, $a_i=p^i_{1i}$, and $b_i=p^i_{1i+1}$. 

For each integer $i$ ($0\leq i\leq D$), the $i$th {\it distance matrix} $A_i$ 
of $\Gamma$ has rows and columns indexed by the vertex of $\Gamma$, and, 
for any $x,y\in V(\Gamma)$, 
\begin{equation*}
(A_i)_{x,y} = \left \{ \begin{aligned}
1\text{~if~}d(x,y)=i,\\
0\text{~if~}d(x,y)\ne i.
\end{aligned}\right.
\end{equation*}

Then $A:=A_1$ is just the {\it adjacency matrix} of $\Gamma$, 
$A_0=I$, $A_i^{\top}=A_i$ ($0\leq i\leq D$), and 
\begin{equation*}
A_iA_j=\sum_{h=0}^D p^h_{ij}A_h ~~~ (0\leq i,j\leq D),
\end{equation*}
in particular,
\begin{equation*}
A_1A_i=b_{i-1}A_{i-1}+a_{i}A_{i}+c_{i+1}A_{i+1} ~~~~ (1\leq i\leq D-1),
\end{equation*}
\begin{equation*}
A_1A_D=b_{D-1}A_{D-1}+a_{D}A_{D}, 
\end{equation*}
and this implies that $A_i=p_i(A_1)$ for certain polynomial $p_i$ of degree $i$. 

The {\it Bose-Mesner} algebra ${\cal M}$ of $\Gamma$ is a matrix algebra generated 
by $A_1$ over ${\mathbb C}$. It follows that ${\cal M}$ has dimension $D+1$, 
and it is spanned by the set of matrices $A_0=I,A_1,\ldots,A_D$, which form a basis of ${\cal M}$.

Since the algebra ${\cal M}$ is semi-simple and commutative, ${\cal M}$ also has 
a basis of pairwise orthogonal idempotents $E_0:=\frac{1}{|V(\Gamma)|}J,E_1,\ldots,E_D$ 
(the so-called {\it primitive idempotents} of ${\cal M}$):
\begin{equation*}
E_iE_j=\delta_{ij}E_i~~(0\leq i,j\leq D),~~E_i=E_i^{\top}~~(0\leq i,j\leq D), 
\end{equation*}
\begin{equation*}
E_0+E_1+\ldots+E_D=I, 
\end{equation*}
where $J$ is the all ones matrix.
 
In fact, $E_j$ ($0\leq j\leq D$) is the matrix representing orthogonal projection onto 
the eigenspace of $A_1$ corresponding to some eigenvalue of $\Gamma$.
In other words, one can write
\begin{equation*}
A_1=\sum_{j=0}^D \theta_jE_j,
\end{equation*} 
where $\theta_j$ ($0\leq j\leq D$) are the real and pairwise distinct 
scalars, known as the {\it eigenvalues} of $\Gamma$.
We say that the eigenvalues are in {\it natural} order if $b_0=\theta_0>\theta_1>\ldots>\theta_D$.
We denote $\hat{\theta}_i=-1-\frac{b_1}{\theta_i+1}$ for $i\in \{1,D\}$. 

The Bose-Mesner algebra ${\cal M}$ is also closed under entrywise (Hadamard or Schur) 
matrix multiplication, denoted by $\circ$. Then the matrices $A_0$, $A_1$, $\ldots$, $A_D$ 
are the primitive idempotents of ${\cal M}$ with respect to $\circ$, i.e., 
$A_i\circ A_j= \delta_{ij}A_i$, and $\sum_{i=0}^D A_i=J$. 
This implies that
\[
E_i\circ E_j=\sum_{h=0}^{D} q_{ij}^h E_h ~~~ (0\leq i,j\leq D)
\]
holds for some real numbers $q_{ij}^h$, known as the {\it Krein parameters} of $\Gamma$. 

Let $\Gamma$ be a distance-regular graph, and $E$ be a primitive idempotent of its Bose-Mesner algebra.
The graph $\Gamma$ is called {\it $Q$-polynomial} (with respect to $E$) 
if there exist real numbers $c_i^*$, $a_i^*$, $b_{i-1}^*$ ($1\leq i\leq D$) 
and an ordering of primitive idempotents 
such that $E_0=\frac{1}{|V(\Gamma)|}J$ and $E_1=E$, and
\[
E_1\circ E_i=b_{i-1}^*E_{i-1} + a_i^*E_i + c_{i+1}^*E_{i+1} ~~~ (1\leq i\leq D-1),
\]
\[
E_1\circ E_D=b_{D-1}^*E_{D-1} + a_D^*E_D.
\]

Note that a $Q$-polynomial ordering of the eigenvalues/idempotents does not have to be the natural ordering.

Further, the {\it dual eigenvalues} of $\Gamma$ {\it associated with} $E$ 
are the real scalars $\theta_i^*$ ($0\leq i\leq D$) defined by 
\begin{equation*}
E=\frac{1}{|V(\Gamma)|}\sum_{i=0}^D \theta_i^* A_i.
\end{equation*} 

Let ${\mathbb V}={\mathbb R}^{V(\Gamma)}$ denote a vector space of columns, 
whose coordinates are indexed by the set $V(\Gamma)$. 
For each vertex $x\in V(\Gamma)$, 
define $\hat{x}\in {\mathbb V}$ by
\begin{equation*}
\hat{x}=(0,0,\ldots,0,1,0,\ldots,0)^{\top},
\end{equation*}
where the 1 is in coordinate $x$. Let $\langle,\rangle$ 
denote the dot product
\begin{equation*}
\langle x,y \rangle = x^{\top}y~~~(x,y\in {\mathbb V}).
\end{equation*}

Fix any vertex $x\in V(\Gamma)$. We will call $x$ a {\it base} vertex.
For each integer $i$ ($0\leq i\leq D$), 
let $E_i^*=E_i^*(x)$ 
denote the diagonal matrix with rows and columns indexed 
by $V(\Gamma)$, and defined by 
\[
(E_i^*)_{y,y}=(A_i)_{x,y}=\langle A_i \hat{x},\hat{y}\rangle ~~~ (y\in V(\Gamma)).
\]

Note that $E_i^*E_j^*=\delta_{ij}E_i^*$, $\sum_{i=0}^DE_i^*=I$.
These matrices span the {\it dual Bose-Mesner} algebra ${\cal M}^* = {\cal M}^*(x)$ 
with respect to $x$. The {\it Terwilliger} (or {\it subconstituent}) algebra 
${\cal T}={\cal T}(x)$ with respect to $x$ is the matrix algebra generated 
by ${\cal M}$ and ${\cal M}^*(x)$, see \cite{SubAlgPaper}. 
(Note that ${\cal T}(x)$ may depend on the base vertex $x$, as for example 
it is the case for the twisted Grassmann graphs, see \cite{TwistG}.)
That this algebra is semi-simple is a standard observation. 
 
A ${\cal T}$-module is any subspace $W\subset {\mathbb V}$ such that $Tw\in W$ for any $T\in {\cal T}$, $w\in W$.
A non-trivial ${\cal T}$-module is {\it irreducible} if it does not properly contain a non-zero 
${\cal T}$-module. Since ${\cal T}$ is semi-simple, each ${\cal T}$-module is a direct sum of 
irreducible ${\cal T}$-modules, and $\mathbb{V}$ decomposes into an orthogonal direct sum of 
irreducible ${\cal T}$-modules. 

Let $W$ be an irreducible ${\cal T}$-module. We define the {\it endpoint} 
of $W$ by ${\rm min}\{i:~E_i^*W\ne 0\}$. An irreducible ${\cal T}$-module $W$ is called {\it thin} 
if ${\rm dim}(E_i^*W)\leq 1$ for all $i=0,1,\ldots,D$; 
the graph $\Gamma$ is called {\it thin} if, for any of its vertices $x$, 
each irreducible ${\cal T}(x)$-module is thin.  
There is a unique irreducible ${\cal T}$-module of endpoint $0$, called 
the {\it trivial} module; it is thin and has basis $\{E_i^*{\bf 1}\mid 0\leq i\leq D\}$, 
where ${\bf 1}$ is the all ones vector.
\medskip

Let us recall the definitions of some families of $Q$-polynomial distance-regular graphs.
The {\it Johnson graph} $J(N,D)$ ($N\geq 2D$) has as vertices all $D$-element subsets of an $N$-element set, 
with two vertices adjacent if their intersection contains $D-1$ elements. 
The {\it binary Hamming graph} $H(N,2)$ (the $N$-cube) has as vertices all words of length $N$ over 
$\{0,1\}$, with two words adjacent if they differ in exactly one position.
The {\it halved graph} of the $N$-cube, denoted by $\frac{1}{2}H(N,2)$, has vertex set consisting 
of all even-weight words of length $N$ over $\{0,1\}$, with two words adjacent if they differ in exactly 
two positions. 
For $\Gamma=\frac{1}{2}H(N,2)$ (with $N$ even) or $\Gamma=J(2D,D)$, the {\it antipodal quotient} 
(the {\it folded graph}) $\widetilde{\Gamma}$ of $\Gamma$ has vertex set consisting of all subsets 
$\{x,y\}$, $x,y\in V(\Gamma)$, where $d(x,y)$ is the diameter of $\Gamma$, with 
$\{x,y\}$ and $\{x',y'\}$ adjacent in $\widetilde{\Gamma}$ if $x\sim_{\Gamma} x'$ 
or $x\sim_{\Gamma} y'$.

We say that a distance-regular graph $\Gamma$ has {\it classical parameters} $(D,b,\alpha,\beta)$ 
if the diameter of $\Gamma$ is $D$, and the intersection numbers of $\Gamma$ 
satisfy 
\begin{equation}\label{classparamc_i}
c_i=\genfrac{[}{]}{0pt}{}{i}{1}\Big(1+\alpha\genfrac{[}{]}{0pt}{}{i-1}{1}\Big),
\end{equation}
so that, in particular, $c_2=(b+1)(\alpha+1)$, 
\begin{equation}\label{classparamb_i}
b_i=\Big(\genfrac{[}{]}{0pt}{}{D}{1}-\genfrac{[}{]}{0pt}{}{i}{1}\Big)\Big(\beta-\alpha\genfrac{[}{]}{0pt}{}{i}{1}\Big),
\end{equation}
where 
$$\genfrac{[}{]}{0pt}{}{j}{1}:=1+b+b^2+\ldots+b^{j-1}.$$

Note that $\Gamma$ is $Q$-polynomial, see \cite[Corollary 8.4.2]{BCN}, and let $E$ 
be the corresponding primitive idempotent of $\Gamma$. The following relation 
(\cite[Eq. (52)]{KiteFree}) between the dual eigenvalues associated with $E$ will be useful:
\begin{equation}\label{classparamtheta_i}
\theta_i^*-\theta_0^* = (\theta_1^*-\theta_0^*)\genfrac{[}{]}{0pt}{}{i}{1}b^{1-i}.
\end{equation}

The Johnson graphs and the Hamming graphs have classical parameters, see \cite[Chapter~6]{BCN}.
It is known (see \cite[Prop.~6.3.1]{BCN}) when a graph with classical parameters $(D,b,\alpha,\beta)$
with diameter $D\geq 3$ is imprimitive (bipartite or antipodal). It is bipartite if and only if $\alpha=0$ and $\beta=1$; 
and it is antipodal if and only if $b=1$ and $\beta=1 + \alpha(D-1)$, in which case it is 
an antipodal double cover of its folded graph. This folded graph has diameter $D'$ 
and intersection numbers defined by
$$b_i=(D-i)(1 + \alpha(D-1-i)),~~c_i=i(1 + \alpha(i-1))\text{~~for~~}i<D',$$
and $c_{D'}=\gamma D'(1 + \alpha(D'-1))$, where $\gamma=1$ if $D=2D'+1$ and
$\gamma=2$ if $D=2D'$. The graphs with such intersection numbers are called {\it pseudo partition} graphs. 
It is also known (see \cite[Theorem~3.3]{BussNeu}) that pseudo partition graphs
with diameter $D'\geq 3$ must have the same intersection arrays as one of the three
following families of partition graphs: the folded cubes ($\alpha=0$), 
the folded Johnson graphs ($\alpha=1$), and the folded halved cubes ($\alpha=2$). 

The folded cubes are determined by their intersection arrays, except for the folded 6-cube, 
see \cite[Prop.~9.2.7]{BCN}. The characterization of the other two families 
of partition graphs has been reduced in \cite{BussNeu}, 
\cite{Metsch97}, \cite{Metsch971}, \cite{Metsch03}
(for the detailed background see Section \ref{SectionType2})
to the problem about the folded halved cubes of diameter 3 and 4 with intersection arrays 
$\{91,66,45;1,6,15\}$, $\{66,45,28;1,6,30\}$, and $\{120,91,66,45;1,6,15,56\}$ respectively.
In this paper, we complete the characterization of pseudo partition graphs 
(see Theorem \ref{MainCoro2} in Section \ref{SectionType2}).
Our proof also shows the previous results by Bussemaker and Neumaier \cite{BussNeu}, and Metsch 
\cite{Metsch97}, \cite{Metsch971}, \cite{Metsch03}.

\section{Triple intersection numbers $[i,i+1,i+1]$}\label{spearSection}

In this section, we show that (\ref{spear1}) holds, see Theorem \ref{spearTheo} below.

Let $\Gamma$ denote a distance-regular graph with diameter $D\geq 3$. 
For all $x,y\in V(\Gamma)$, and all integers $i$ and $j$, define 
$$p_{ij}(x,y):=\sum_{z\in \Gamma_i(x)\cap \Gamma_j(y)}\hat{z},$$
and 
$$x_y^-:=p_{1h-1}(x,y),~x_y^0:=p_{1h}(x,y),~x_y^+:=p_{1h+1}(x,y),$$
where $h=d(x,y)$.

The following lemma is just Lemma 2.9 from \cite{KiteFree} written 
in a slightly different form.

\begin{lemma}\label{Lemma2.9KiteFree} 
Let $\Gamma$ denote a distance-regular graph 
with diameter $D\ge 3$, suppose that $\Gamma$ is $Q$-polynomial 
with respect to the primitive idempotent 
$$E_1=|V(\Gamma)|^{-1}\sum_{h=0}^D \theta_{h}^*A_h.$$

Then, for any adjacent vertices $x,y\in V(\Gamma)$, and any integer 
$i$ ($1\le i\le D$), the vector 
\[
\frac{1}{p_{ii+1}^1} p_{ii+1}(x,y) 
+ \hat{x}\frac{\theta_1^*-\theta_{i}^*}{\theta_0^*-\theta_1^*} 
- \hat{y}\frac{(\theta_1^*-\theta_2^*)(\theta_1^*-\theta_{i}^*)-(\theta_0^*-\theta_1^*)(\theta_2^*-\theta_{i+1}^*)}{(\theta_0^*-\theta_1^*)(\theta_0^*-\theta_2^*)}
- \frac{1}{b_1}x_{y}^{+}\frac{\theta_0^*+\theta_1^*-\theta_i^*-\theta_{i+1}^*}{\theta_0^*-\theta_2^*}
\]
is orthogonal to $E_0{\mathbb V}+E_1{\mathbb V}$.
\end{lemma}
\proof By \cite[Lemma 2.9]{KiteFree}, we have that the vector
$$
\frac{1}{p_{ii-1}^1} p_{ii-1}(x,y) 
- \hat{x}\frac{(\theta_1^*-\theta_2^*)(\theta_1^*-\theta_{i-1}^*)-(\theta_0^*-\theta_1^*)(\theta_2^*-\theta_{i}^*)}{(\theta_0^*-\theta_1^*)(\theta_0^*-\theta_2^*)}
+ \hat{y}\frac{\theta_1^*-\theta_{i-1}^*}{\theta_0^*-\theta_1^*}
- \frac{1}{b_1}y_{x}^{+}\frac{\theta_0^*+\theta_1^*-\theta_{i-1}^*-\theta_{i}^*}{\theta_0^*-\theta_2^*}
$$
is orthogonal to $E_0{\mathbb V}+E_1{\mathbb V}$.

By the definition of $p_{ij}(x,y)$, we have $p_{ii-1}(x,y)=p_{i-1i}(y,x)$.
Now the lemma follows by changing the roles of $x$ and $y$ and replacing $i$ by $i+1$.\wbull









\medskip

The following theorem is an analogue of Theorem 2.11 from \cite{KiteFree}.

\begin{theo}\label{spearTheo} 
Let $\Gamma$ denote a distance-regular graph 
with diameter $D\ge 3$, suppose that $\Gamma$ is $Q$-polynomial 
with respect to the primitive idempotent 
$$E_1=|V(\Gamma)|^{-1}\sum_{h=0}^D \theta_{h}^*A_h.$$


Then the following holds.

\quad $(1)$ If $y\sim z$ then 
$$[i,i+1,i+1] = p_{i,i+1}^1\Big(\frac{[1,2,2]}{b_1}\times \frac{(\theta_{2}^*-\theta_{1}^*)(\theta_0^*+\theta_1^*-\theta_i^*-\theta_{i+1}^*)}{(\theta_0^*-\theta_2^*)(\theta_{i+1}^*-\theta_{i}^*)} + \frac{\theta_1^*-\theta_{i}^*}{\theta_{i+1}^*-\theta_{i}^*}\Big).$$

\quad $(2)$ If $y\not\sim z$ then
\begin{multline}\label{spearnotsim}
[i,i+1,i+1] = p_{i,i+1}^1\Big(\frac{[1,2,2]}{b_1}\times 
\frac{(\theta_2^*-\theta_1^*)(\theta_0^*+\theta_1^*-\theta_i^*-\theta_{i+1}^*)}{(\theta_0^*-\theta_2^*)(\theta_{i+1}^*-\theta_i^*)} + 
\frac{(\theta_1^*-\theta_i^*)}{(\theta_{i+1}^*-\theta_i^*)} + {} \\
{} + \frac{(\theta_{0}^*-\theta_1^*)(\theta_{1}^*-\theta_2^*)(\theta_{2}^*-\theta_{i+1}^*)- (\theta_{1}^*-\theta_2^*)^2(\theta_{1}^*-\theta_i^*)}{(\theta_{0}^*-\theta_1^*)(\theta_{0}^*-\theta_2^*)(\theta_{i+1}^*-\theta_i^*)} + {} \\ 
{} + \frac{1}{b_1}\times \frac{(\theta_{0}^*-\theta_1^*)(\theta_0^*+\theta_1^*-\theta_i^*-\theta_{i+1}^*) - c_i(\theta_0^*-\theta_2^*)(\theta_{i-1}^*-\theta_i^*)}{(\theta_0^*-\theta_2^*)(\theta_{i+1}^*-\theta_i^*)}\Big).
\end{multline}
\end{theo} 
\proof Denote $|V(\Gamma)|$ by $v$ for short.
We follow the proof of \cite[Theorem 2.11]{KiteFree}. 
Let us compute the inner product of $E_1\hat{z}$ and the vector 
in Lemma \ref{Lemma2.9KiteFree}, and set the result equal to 0. 
Note that 
$$v\langle E_1\hat{z},p_{ii+1}(x,y)\rangle=
[i,i+1,i-1]\theta_{i-1}^*+[i,i+1,i]\theta_{i}^*+[i,i+1,i+1]\theta_{i+1}^*,$$
and $[i,i+1,i-1]+[i,i+1,i]+[i,i+1,i+1]=p_{i,i+1}^1$.

$(1)$ If $y\sim z$ then $[i,i+1,i-1]=0$, and 
$$v\langle E_1\hat{z},p_{ii+1}(x,y)\rangle=
(p_{i,i+1}^1-[i,i+1,i+1])\theta_{i}^*+[i,i+1,i+1]\theta_{i+1}^*.$$

Further, 
$$v\langle E_1\hat{z},\hat{x}\rangle=\theta_1^*,$$
$$v\langle E_1\hat{z},\hat{y}\rangle=\theta_1^*,$$
$$v\langle E_1\hat{z},x_{y}^{+}\rangle=
b_1\theta_{1}^*+[1,2,2](\theta_{2}^*-\theta_{1}^*).$$

It now follows that
$$[i,i+1,i+1] = p_{i,i+1}^1\Big(\frac{[1,2,2]}{b_1}\times \frac{(\theta_{2}^*-\theta_{1}^*)(\theta_0^*+\theta_1^*-\theta_i^*-\theta_{i+1}^*)}{(\theta_0^*-\theta_2^*)(\theta_{i+1}^*-\theta_{i}^*)} + \frac{\theta_1^*-\theta_{i}^*}{\theta_{i+1}^*-\theta_{i}^*}\Big).$$

$(2)$ If $y\not\sim z$ then $[i,i+1,i-1]=p_{i+1,i-1}^{2}$, and 
$$v\langle E_1\hat{z},p_{ii+1}(x,y)\rangle=
p_{i+1,i-1}^2\theta_{i-1}^*+(p_{i,i+1}^1-p_{i+1,i-1}^2-[i,i+1,i+1])\theta_{i}^*+[i,i+1,i+1]\theta_{i+1}^*.$$

Further, 
$$v\langle E_1\hat{z},\hat{x}\rangle=\theta_1^*,$$
$$v\langle E_1\hat{z},\hat{y}\rangle=\theta_2^*,$$
$$v\langle E_1\hat{z},x_{y}^{+}\rangle=[1,2,0]\theta_{0}^*+[1,2,1]\theta_{1}^*+[1,2,2]\theta_{2}^*=\theta_{0}^*+(b_1-1-[1,2,2])\theta_{1}^*+[1,2,2]\theta_{2}^*,$$
which gives (\ref{spearnotsim}). The theorem is proved.\wbull



\medskip

Using (\ref{classparamc_i})--(\ref{classparamtheta_i}), we obtain:

\begin{coro}
Let $\Gamma$ denote a distance-regular graph 
with diameter $D\ge 3$. Suppose that $\Gamma$ is $Q$-polynomial 
with respect to the primitive idempotent 
$$E_1=|V(\Gamma)|^{-1}\sum_{h=0}^d \theta_{h}^*A_h$$
and has classical parameters $(D,b,\alpha,\beta)$.

Then the following holds.

\quad $(1)$ If $y\sim z$ then 
$$[i,i+1,i+1] = p_{i,i+1}^1\Big(\frac{[1,2,2]}{b_1}\times \frac{2\genfrac{[}{]}{0pt}{}{i+1}{1}-(1+b^i)}{1+b} - b\genfrac{[}{]}{0pt}{}{i-1}{1}\Big).$$

\quad $(2)$ If $y\not\sim z$ then
\begin{multline*}
[i,i+1,i+1] = p_{i,i+1}^1\Big(\frac{[1,2,2]}{b_1}\times \frac{2\genfrac{[}{]}{0pt}{}{i+1}{1}-(1+b^i)}{1+b} - b\genfrac{[}{]}{0pt}{}{i-1}{1} + \frac{b}{b_1}\times 
\big(c_i - \frac{(2\genfrac{[}{]}{0pt}{}{i+1}{1}-(1+b^i))}{1+b}\big) \Big).
\end{multline*}
\end{coro}

\section{The Terwilliger polynomial}

Let $\Gamma$ be a distance-regular graph with diameter $D\geq 3$, 
suppose that $\Gamma$ is $Q$-polynomial with $Q$-polynomial ordering $E_0,E_1,\ldots,E_D$. 
Fix a base vertex $x\in V(\Gamma)$, ${\cal T}={\cal T}(x)$ with respect to $x$, 
$E_i^*=E_i^*(x)$ and denote $\widetilde{A}_i=E_i^*A_1E_i^*$ for  
$i=0,1,\ldots,D$. We write $\widetilde{A}$ instead of $\widetilde{A}_1$ for short. 
For notational convenience, we also set $\widetilde{A}^0:=E_1^*$ and $\widetilde{J}:=E_1^*JE_1^*$.
With an appropriate ordering of the vertices of $\Gamma$, one can see 
that 
$$
\widetilde{A}=\left(
 \begin{tabular}{ll}
 $N$ & 0\\
 0 & 0
 \end{tabular}
 \right), 
$$
where the principal submatrix $N$ is, in fact, the adjacency matrix of $\Gamma(x)$, the local graph of $x$.

We now recall some facts about irreducible ${\cal T}$-modules of endpoint 1, see \cite{SN}, cf. \cite{IH}.
Let $U_1^*$ be the the subspace of $E_1^*\mathbb{V}$, which is orthogonal to ${\bf 1}$.
Let $W$ be an irreducible ${\cal T}(x)$-module of endpoint $1$. Then $E_1^*W$ 
is a one-dimensional subspace of $U_1^*$; in particular, any non-zero vector $w\in E_1^*W$ 
is an eigenvector of $\widetilde{A}$, and $W={\cal T}w$. Conversely, for an eigenvector 
$w$ of $\widetilde{A}$, the subspace $W={\cal T}w$ is an irreducible ${\cal T}$-module of endpoint 1.
Let $a_0(W)$ denote the corresponding eigenvalue. Note that $a_0(W)$ is a non-principal eigenvalue 
of the local graph of $x$.

Set
$$w_i^+=E_i^*A_{i-1}E_1^*w,~~w_i^-=E_i^*A_{i+1}E_1^*w~~(1\leq i\leq D).$$

It follows from \cite{SN} that $W$ is thin if and only if $w_i^+,w_i^-$ are linearly dependent for every $i$, 
$2\leq i\leq D-1$. Note that the determinant of a Gram matrix of $w_i^+,w_i^-$ is non-negative:
\begin{equation}\label{DetP}
{\rm det}\left(
 \begin{tabular}{ll}
 $\langle w_i^+,w_i^+\rangle$ & $\langle w_i^+,w_i^-\rangle$\\
 $\langle w_i^+,w_i^-\rangle$ & $\langle w_i^-,w_i^-\rangle$
 \end{tabular}
 \right)\geq 0.
\end{equation}

Terwilliger \cite{SN} implicitly showed that there exist polynomials 
$p_i^{\epsilon\delta}\in {\mathbb C}[\lambda]$, $\epsilon,\delta\in\{+,-\}$, 
which depend only on $i$ and the intersection numbers of $\Gamma$,  
such that 
$$
E_1^*A_{i\epsilon 1}E_i^*A_{i\delta 1}E_1^* = \alpha_{\epsilon\delta} {\widetilde{J}}+p_i^{\epsilon\delta}(\widetilde{A})
$$
for some $\alpha_{\epsilon\delta}\in {\mathbb C}$, and therefore
$$\langle w_i^{\epsilon},w_i^{\delta}\rangle = \|w\|^2 p_i^{\epsilon\delta}(a_0(W)),$$
as 
$$
\langle w_i^{\epsilon},w_i^{\delta}\rangle=w^{\top}E_1^*A_{i\epsilon 1}E_i^*A_{i\delta 1}E_1^*w,
\text{~~and~~}w^{\top}\widetilde{J}w=0.
$$

Define 
\begin{equation}\label{TerwPoly}
T_i(\lambda):=p_i^{++}(\lambda)p_i^{--}(\lambda)-p_i^{+-}(\lambda)^2.
\end{equation}

Taking into account (\ref{DetP}), one has the following result.

\begin{theo}\label{TerwPolyTheo}
Let $\Gamma$ be a $Q$-polynomial distance-regular graph with diameter $D\geq 3$. 
Then, for any $i=2,\ldots,D-1$, for any vertex $x\in V(\Gamma)$ and any non-principal eigenvalue $\eta$ 
of the local graph of $x$, $T_i(\eta)\geq 0$ holds, with equality if and only if 
${\cal T}(x)w$ is a thin module of endpoint $1$, where $w$ is an eigenvector of $\widetilde{A}$ 
with eigenvalue $\eta$.
\end{theo}

We will call the polynomial $T_i(\lambda)$ the {\it Terwilliger polynomial} of $\Gamma$.

\begin{remark}\label{Tclaimremark}
\begin{itemize}
\item[(i)] Terwilliger \cite{SN} claimed that 
\begin{equation}\label{Tclaim}
T_i(\lambda)\text{~is~independent~on~}i\text{~up~to~a~scalar~multiple~}.
\end{equation}
However, the proof was not given, and it seems 
to be non-trivial. Therefore, the independence $T_i(\lambda)$ on $i$ should 
be explained/verified (although we believe that this is true, and this 
was true for all the examples we calculated the polynomials for).  
See also Remark \ref{Tclaimdiscuss} below.

\item[(ii)] One might expect that $p_i^{++},p_i^{--},p_i^{+-}$ do not depend 
on $i$ up to a scalar multiple. This is not true in general.

\item[(iii)] In this section, we calculate $T_i(\lambda)$ explicitly 
for the first possible value of $i$, i.e. $i=2$. Therefore, in what follows, 
we omit the sub-index $i$ and write just $T$, $p^{++}$, $p^{--}$, $p^{+-}$ 
(assuming that $i=2$).
\end{itemize}
\end{remark}

\begin{remark}\label{Tclaimdiscuss}
Suppose that $\Gamma$ has classical parameters. 
We guess that in this case one may use the following approach to prove (\ref{Tclaim}).
Recall the definitions of the {\it raise} operator $R$:
\[R:=\sum_{i=0}^{D-1}E_{i+1}^*A_1E_{i}^*,\]
and the {\it lower} operator $L$:
\[L:=\sum_{i=0}^{D-1}E_{i}^*A_1E_{i+1}^*,\]
which are well known in the theory of Terwilliger algebras. Note that $L^{\top}=R$.

If $W$ is a non-thin module with endpoint $1$ then, with the above notation, it follows from 
\cite[Lemma 2.6]{IH} and \cite[Theorem 5.1]{IH} that 
\[Rw_{i}^+=c_iw_{i+1}^+,\]
\[Rw_{i}^-=r_i^+w_{i+1}^++r_i^-w_{i+1}^-,\]
\[Lw_{i+1}^-=b_{i+1}w_{i}^-,\]
\[Lw_{i+1}^+=l_{i+1}^+w_{i}^++l_{i+1}^-w_{i}^-,\]
for some $r_i^{\pm}, l_{i+1}^{\pm}$, which can be considered as 
polynomials in $\eta:=a_0(W)$ (see \cite{IH} for the details), and $2\leq i\leq D-2$.

Further, using these relations and calculating the inner product $\langle Rw_i^+,Rw_i^-\rangle$ 
in two different ways: 
\begin{equation}\label{2ways}
\langle Rw_i^+,Rw_i^-\rangle=(w_i^+)^{\top}R^{\top}Rw_i^-=(w_i^-)^{\top}R^{\top}Rw_i^+,
\end{equation}
one can get the equality
\begin{equation}\label{2wayseq}
\langle w_i^+,w_i^+\rangle\cdot (r_i^+l_{i+1}^+) + \langle w_i^+,w_i^-\rangle\cdot (r_i^+l_{i+1}^-+r_i^-b_{i+1})=\langle w_i^-,w_i^-\rangle\cdot (c_il_{i+1}^-) + \langle w_i^+,w_i^-\rangle\cdot (c_i^+l_{i+1}^+),
\end{equation}
and then express 
\begin{equation*}
{\rm det}\left(
 \begin{tabular}{ll}
 $\langle w_{i+1}^+,w_{i+1}^+\rangle$ & $\langle w_{i+1}^+,w_{i+1}^-\rangle$\\
 $\langle w_{i+1}^+,w_{i+1}^-\rangle$ & $\langle w_{i+1}^-,w_{i+1}^-\rangle$
 \end{tabular}
 \right) 
 \text{\rm~~~in~terms~of~~~}
 {\rm det}\left(
 \begin{tabular}{ll}
 $\langle w_i^+,w_i^+\rangle$ & $\langle w_i^+,w_i^-\rangle$\\
 $\langle w_i^+,w_i^-\rangle$ & $\langle w_i^-,w_i^-\rangle$
 \end{tabular}
 \right),
\end{equation*}
and a constant multiple. But the details still need to be checked.

Below we will make use of the results from Section \ref{spearSection} to calculate 
polynomials $T_2(\lambda)$ and $p_2^{--}(\lambda)$. 

Note that (\ref{2wayseq}) allows us to calculate $\langle w_i^-,w_i^-\rangle$ 
as a polynomial in $\eta$, if $\langle w_i^+,w_i^-\rangle=\|w\|^2p_i^{+-}(\eta)$ 
and $\langle w_i^+,w_i^+\rangle=\|w\|^2p_i^{++}(\eta)$ are known. 
(The polynomials $p_2^{++}(\lambda)$ and $p_2^{+-}(\lambda)$ 
can be relatively easily determined.)

We do not know whether for any distance-regular graph with classical parameters 
the polynomial $p_2^{--}(\lambda)$ calculated using Equation (\ref{2wayseq}) is the same 
as the one in Lemma \ref{pmimi}. But if they are not the same this would give strong new information. 
The above approach probably will also work for all $Q$-polynomial distance-regular graphs. 
\end{remark}

The polynomial $p_i^{+-}$ can be easily determined for any $i$.

\begin{lemma}\label{pplmi}
The following holds:
$$E_1^*A_{i-1}E_i^*A_{i+1}E_1^*=(\widetilde{J}-\widetilde{A}-E_1^*)p^2_{i-1,i+1},~~(2\leq i\leq D-1)$$
in particular, $p_i^{+-}(\lambda)=-p^2_{i-1,i+1}(\lambda+1)$.
\end{lemma}
\proof 
For a pair of vertices $y,z\in \Gamma(x)$, we have 
$(E_1^*A_{i-1}E_i^*A_{i+1}E_1^*)_{y,z} = |\Gamma_i(x)\cap \Gamma_{i-1}(y)\cap \Gamma_{i+1}(z)|$, i.e., 
\begin{equation}\label{eq_for_pplusplus2}
(E_1^*A_{i-1}E_i^*A_{i+1}E_1^*)_{y,z} =\left \{ \begin{aligned}
0\text{~if~}y=z,\\
0\text{~if~}y\sim z,\\
p_{i-1,i+1}^2\text{~if~}y\not\sim z,~y\ne z,
\end{aligned}\right.
\end{equation}
which shows the lemma.\wbull

\begin{lemma}\label{pplpl}
The following holds:
\begin{equation}\label{eq_for_pplusplus}
E_1^*AE_2^*AE_1^* = (c_2-1)\widetilde{J} + (k-c_2)E_1^* + (a_1-c_2)\widetilde{A} - \widetilde{A}^2,
\end{equation}
in particular, $p^{++}(\lambda)=-\lambda^2 + (a_1-c_2)\lambda + (k-c_2)$. 
\end{lemma}
\proof For a pair of vertices $y,z\in \Gamma(x)$, we have 
$(E_1^*AE_2^*AE_1^*)_{y,z} = |\Gamma_2(x)\cap \Gamma(y)\cap \Gamma(z)|$, 
and note that $\widetilde{A}_{y,z}^2=|\Gamma(x)\cap \Gamma(y)\cap \Gamma(z)|$
so that 

\begin{equation}\label{eq_for_pplusplus2}
(E_1^*AE_2^*AE_1^*)_{y,z} =\left \{ \begin{aligned}
b_1\text{~if~}y=z,\\
a_1-1-\widetilde{A}_{y,z}^2\text{~if~}y\sim z,\\
c_2-1-\widetilde{A}_{y,z}^2\text{~if~}y\not\sim z,~y\ne z.
\end{aligned}\right.
\end{equation}
 
The lemma follows by comparing (\ref{eq_for_pplusplus2}) with 
the right-hand side of (\ref{eq_for_pplusplus}).\wbull

\begin{lemma}\label{pmimi} 
The following holds:
\begin{equation}\label{eq_for_pminusminus}
E_1^*A_3E_2^*A_3E_1^* = p^1_{23}\Big(\tau_2\widetilde{J} + (1-a_1\tau_0-\tau_2)E_1^* 
+ (\tau_1-\tau_2)\widetilde{A} + \tau_0\widetilde{A}^2\Big),
\end{equation}
where 
$$\tau_0=\frac{1}{b_1}\times \frac{(\theta_{2}^*-\theta_{1}^*)(\theta_0^*+\theta_1^*-\theta_2^*-\theta_{3}^*)}{(\theta_0^*-\theta_2^*)(\theta_{3}^*-\theta_{2}^*)},$$
$$\tau_1=\frac{\theta_2^*-\theta_{1}^*}{\theta_{3}^*-\theta_{2}^*}\Big( 
\frac{(\theta_1^*-\theta_{3}^*)}{(\theta_0^*-\theta_2^*)} - \frac{a_1-1}{b_1} - 
\frac{(a_1-1)(\theta_1^*-\theta_{3}^*)}{b_1(\theta_0^*-\theta_2^*)}\Big),$$
\begin{multline}
\tau_2=\frac{\theta_2^*-\theta_{1}^*}{\theta_{3}^*-\theta_{2}^*}\Big( 
\frac{(\theta_1^*-\theta_{3}^*)}{(\theta_0^*-\theta_2^*)} - \frac{a_1+1-c_2}{b_1} - 
\frac{(a_1+1)(\theta_1^*-\theta_{3}^*)}{b_1(\theta_0^*-\theta_2^*)} + {} \\
{} + \frac{(\theta_{1}^*-\theta_2^*)^2-(\theta_{0}^*-\theta_1^*)(\theta_{2}^*-\theta_{3}^*)}{(\theta_{0}^*-\theta_1^*)(\theta_{0}^*-\theta_2^*)} -
\frac{1}{b_1}\times \frac{(\theta_{0}^*-\theta_1^*)(\theta_0^*+\theta_1^*-\theta_2^*-\theta_{3}^*)}{(\theta_0^*-\theta_2^*)(\theta_{1}^*-\theta_2^*)}\Big)
\end{multline}

In particular, $p^{--}(\lambda)=p^1_{23}\Big(\tau_0\lambda^2 + (\tau_1-\tau_2)\lambda + 
(1-a_1\tau_0-\tau_2)\Big)$. 
\end{lemma}
\proof For a pair of vertices $y,z\in \Gamma(x)$, we have 
$(E_1^*A_3E_2^*A_3E_1^*)_{y,z} = |\Gamma_2(x)\cap \Gamma_3(y)\cap \Gamma_3(z)|=[2,3,3]$. 
Recall that $\widetilde{A}_{y,z}^2=[1,1,1]$ in our notation.

We note that if $y\sim z$ then 
$[1,1,1]+[1,1,2]=a_1-1$ and $[1,1,2]+[1,2,2]=b_1$
hold. This gives 
\begin{equation}\label{1221111}
[1,2,2]=[1,1,1]+b_1-a_1+1\text{~~~if~}y\sim z.
\end{equation} 

If $y\not\sim z$ and $y\ne z$ then 
$[1,1,1]+[1,1,2]=a_1$ and $[1,1,2]+[1,2,2]=b_1-1$ 
so that 
\begin{equation}\label{1221112}
[1,2,2]=[1,1,1]+b_1-a_1-1\text{~~~if~}y\not\sim z,~y\ne z.
\end{equation} 

By Theorem \ref{spearTheo} and (\ref{1221111}), (\ref{1221112}), 
we see that 
\begin{equation}\label{eq_for_pminusminus2}
\frac{[2,3,3]}{p^1_{23}} = \left \{ \begin{aligned}
1\text{~if~}y=z,\\
\tau_0\widetilde{A}_{y,z}^2 + \tau_1\text{~if~}y\sim z,\\
\tau_0\widetilde{A}_{y,z}^2 + \tau_2\text{~if~}y\not\sim z,~y\ne z,
\end{aligned}\right.
\end{equation}
where $\tau_0$, $\tau_1$, and $\tau_2$ are rational expressions 
in the dual eigenvalues and intersection numbers of $\Gamma$.
Using (\ref{1221111}) and (\ref{1221112}), we obtain explicit 
expressions for $\tau_0$, $\tau_1$, and $\tau_2$.

The lemma follows by comparing (\ref{eq_for_pminusminus2}) with 
the right-hand side of (\ref{eq_for_pminusminus}).\wbull

\begin{lemma}\label{Tpolyclassical}
Suppose that $\Gamma$ has classical parameters $(D,b,\alpha,\beta)$. 
Then the Terwilliger polynomial of $\Gamma$ is 
\begin{multline}\label{TpolyclassicalEq}
T(\lambda)=\frac{b_2}{\alpha+1}\Big(-\lambda^2 + \lambda\big(\alpha\genfrac{[}{]}{0pt}{}{D}{1}+\beta-\alpha-1-(\alpha+1)(b+1)\big) + \beta\genfrac{[}{]}{0pt}{}{D}{1}-(\alpha+1)(b+1)\Big)\times {} \\
{} \times \Big(\lambda^2+\lambda(2-\alpha b)-\alpha b+1\Big) - b_2^2(\lambda+1)^2.
\end{multline}

Furthermore, the roots of $T(\lambda)$ are
$$\beta-\alpha-1,~-1,~-b-1,~\alpha b\frac{b^{D-1}-1}{b-1}-1.$$
\end{lemma}
\proof Using (\ref{classparamc_i})--(\ref{classparamtheta_i}) and Lemma \ref{pmimi}, 
we obtain the following expressions for $\tau_0$, $\tau_1$, $\tau_2$:
$$\tau_0=\frac{b+1}{b_1},~
\tau_1=\frac{b+1}{b_1}(b_1-a_1+1)-b,~
\tau_2=\frac{b+1}{b_1}(b_1-a_1-1+\alpha b)-b,$$
which gives
\begin{equation}\label{pmimiclass}
p^{--}(\lambda)=p^1_{23}\frac{b+1}{b_1}\big(\lambda^2+\lambda(2-\alpha b)-\alpha b+1\big).
\end{equation}

As $p^1_{23}=\frac{b_1b_2}{c_2}$ and $c_2=(b+1)(\alpha+1)$, we see that 
$p^1_{23}\frac{b+1}{b_1}=\frac{b_2}{\alpha+1}$.

Further, we have $a_1=\alpha\genfrac{[}{]}{0pt}{}{D}{1}+\beta-\alpha-1$, 
$k=\beta\genfrac{[}{]}{0pt}{}{D}{1}$ and $c_2=(\alpha+1)(b+1)$ so that 
\begin{equation}\label{pplplclass}
p^{++}(\lambda)=-\lambda^2 + \lambda\big(\alpha\genfrac{[}{]}{0pt}{}{D}{1}+\beta-\alpha-1-(\alpha+1)(b+1)\big) + \beta\genfrac{[}{]}{0pt}{}{D}{1}-(\alpha+1)(b+1).
\end{equation}

Combining (\ref{pmimiclass}) and (\ref{pplplclass}), we obtain (\ref{TpolyclassicalEq}). 

It is easily seen from (\ref{pmimiclass}) that $p^{--}(-1)=0$ as well as $p^{+-}(-1)=0$.
Therefore, $-1$ is a root of $T(\lambda)$. Substituting $\lambda=\beta-\alpha-1$, 
$-b-1$, or $\alpha b\frac{b^{D-1}-1}{b-1}-1$ into (\ref{TpolyclassicalEq}) 
(after some tedious calculations) shows the lemma.\wbull

\begin{remark}
For any $Q$-polynomial distance-regular graph Terwilliger \cite{SubAlgPaper} calculated 
the possible eigenvalues of the local graph belonging to thin irreducible modules of endpoint $1$. 
We can show that the roots of the Terwilliger polynomial are consistent with his result.
\end{remark}

\section{Pseudo-partition graphs and $Q$-polynomial distance-regular graphs of type 2}\label{SectionType2}

As we mentioned above, the pseudo-partition graphs with diameter at least $3$ 
must have the same intersection arrays as one of the three families of partition graphs: 
the folded cubes, the folded Johnson graphs, and the folded halved cubes. 
These graphs are $Q$-polynomial, and the graphs from the last two families 
are said to be type 2 graphs (see below for the definition).

Recall that a distance-regular graph $\Gamma$ with diameter $D\ge 3$ is said to be 
$Q$-polynomial of type 2 ({\it of type {\rm 2}} for short) if there exist 
$h,h^*,x,y,t^*\in {\mathbb C}$ such that:
\begin{itemize}
\item[(1)] the intersection numbers of $\Gamma$ are given by:
\begin{equation}\label{type2def1}
c_i=\frac{hi(i-t+x)(i-t+y)(i-t+D)}{(2i-t)(2i-t-1)},~1\le i\le D-1,
\end{equation}
\begin{equation}\label{type2defcd}
c_D=\frac{hD(D-t+x)(D-t+y)}{2D-t-1},
\end{equation}
\begin{equation}\label{type2defbi}
b_i=\frac{h(i-t)(i-x)(i-y)(i-D)}{(2i-t)(2i-t+1)},~0\le i\le D-1,
\end{equation}
\begin{equation}\label{type2def2}
b_0=\frac{hxyD}{t-1}.
\end{equation}
\item[(2)] the distinct eigenvalues of $\Gamma$ are given by:
\begin{equation}\label{type2theta}
\theta_i=b_0+hi(i-t^*),~0\le i\le D,
\end{equation}
where 
\begin{equation}\label{ttxyd}
t+t^*=x+y+D+1.
\end{equation}
\item[(3)] the dual intersection numbers $c^*_{i},b^*_{i-1}$ ($1\le i\le D$) 
and dual eigenvalues $\theta^*_i$ ($0\le i\le D$) are obtained by 
replacing $h,t$, and $b_0$ above by $h^*,t^*$, and $b^*_0$, respectively.
\end{itemize}
\medskip

\begin{remark}\label{knowntype2}
Since $c_1$ and $c_1^*$ equal $1$, $h$ and $h^*$ and hence all intersection numbers 
are determined by from $x,y,D$, and $t$. The only known graphs of type $2$ with diameter $D\geq 3$ 
are the following:
\begin{itemize}
\item[(i)] the antipodal quotient of the Johnson graph $J(2t,t)$ 
(the folded Johnson graph) with $D=\frac{t-1}{2}$, $\{x,y\}=\{\frac{t}{2},t\}$ if $t$ is odd, 
and  $D=\frac{t}{2}$, $\{x,y\}=\{\frac{t-1}{2},t\}$ if $t$ is even (in both cases $h=4$),
\item[(ii)] the halved graph $\frac{1}{2}H(2D+1,2)$ of the $(2D+1)$-cube with $t=D+\frac{1}{2}$, 
$\{x,y\}=\{\frac{t-1}{2},\frac{t}{2}\}$, and $h=8$,
\item[(iii)] the antipodal quotient of $\frac{1}{2}H(2t,2)$ (the folded halved cube) 
with $D=\frac{t-1}{2}$, $\{x,y\}=\{\frac{t}{2},t-\frac{1}{2}\}$ if $t$ is odd, 
and $D=\frac{t}{2}$, $\{x,y\}=\{\frac{t-1}{2},t-\frac{1}{2}\}$ if $t$ is even (in both cases $h=8$).
\end{itemize} 

Note that, for the folded Johnson graphs and the folded halved cubes, 
the $Q$-polynomial ordering of primitive idempotents is natural, while the halved $(2D+1)$-cube 
has the $Q$-polynomial ordering $E_0,E_2,E_4,\ldots,E_3,E_1$, 
where $\theta_0>\theta_1>\ldots>\theta_D$ is the natural ordering of the corresponding eigenvalues.
\end{remark}

\begin{remark}
According to \cite[Note~1]{Terw86}, we should note that $t\notin \{1,2,\ldots,2D-1\}$. 
We also mention the correspondence with the notation from the monograph by Bannai and Ito \cite{BI}
(also note that \cite[Note~3]{Terw86} contains misprints: 
there should be $-1-r_i$ instead of $1-r_i$):
$$x=-1-r_1,~y=-1-r_2,~D=-1-r_3,~t=-1-s^*,~t^*=-1-s.$$
\end{remark}

Terwilliger \cite{Terw86} obtained the following result.

\begin{theo}(\cite[Theorem~2.2]{Terw86})\label{Terw86theo}
A type $2$ graph with diameter $D\geq 14$ is either:
\begin{itemize}
\item[(i)] the antipodal quotient of the Johnson graph $J(4D,2D)$ or $J(4D+2,2D+1)$,
\item[(ii)] the halved graph $\frac{1}{2}H(2D+1,2)$ of the $(2D+1)$-cube,
\item[(iii)] the antipodal quotient of $\frac{1}{2}H(4D,2)$ or $\frac{1}{2}H(4D+2,2)$, 
\item[(iv)] a graph not listed above, but with the same intersection array as $(i)$ or $(iii)$.
\end{itemize} 
\end{theo}

In \cite{Neu85}, Neumaier showed that the halved graph $\frac{1}{2}H(2D+1,2)$ of the $(2D+1)$-cube, $D\ge 3$, 
is uniquely determined as distance-regular graph by its intersection array.

In \cite{BussNeu}, Bussemaker and Neumaier proved that the folded Johnson graphs and the folded halved cubes 
with diameter at least 154 are uniquely determined as distance-regular graphs by their intersection arrays. 
Metsch \cite{Metsch97}, \cite{Metsch971}, \cite{Metsch03} improved this result by showing that the same remains 
true for the folded Johnson graphs with diameter at least 3 and the folded halved cubes with diameter at least 5. 
In particular, there are no graphs in (iv) of Theorem \ref{Terw86theo}.

The folded halved cubes of diameter 3 and 4 have intersection arrays $\{91,66,45;1,6,15\}$, 
$\{66,45,28;1,6,30\}$, and $\{120,91,66,45;1,6,15,56\}$ respectively, and the problem of 
their characterization by these intersection arrays remained open. 
In this section, we solve this problem.

\begin{theo}\label{MainCoro2}
The folded Johnson graphs, the folded halved cubes with diameter at least $3$ 
are uniquely determined as distance-regular graphs by their intersection arrays.
\end{theo}

Theorem \ref{MainCoro2} is an immediate consequence of a slightly more general result, 
see Theorem \ref{MainTheo2} below.

Terwilliger \cite{Terw86} noticed that all graphs in Remark \ref{knowntype2} with $D\geq 4$ satisfy 
\begin{equation}\label{graph2condition}
c_3-3c_2+3=b_2-2b_1+k-c_2+2=0,
\end{equation}
however, we should notice that this equation is violated for 
the folded Johnson graph and one of the folded halved cubes with diameter 3.

Following \cite{Terw86}, let us define
\[\gamma_r=\sum_{i=0}^r(-1)^ic_{r-i}{r\choose i},~\beta_r=\sum_{i=0}^r(-1)^ib_{r-i}{r\choose i},\]
so that (\ref{graph2condition}) is equivalent to $\gamma_3=\beta_2-\gamma_2=0$.
Furthermore, if (\ref{graph2condition}) holds, then the intersection numbers of 
the graphs of type 2 take a simple form, see \cite[Corollary~2.5]{Terw86}. 
In particular, it follows from \cite[Corollary~2.5]{Terw86} that 
(\ref{graph2condition}) implies that $\{x,y,D\}=\{\frac{t-1}{2},\frac{t}{2},t-1+\frac{2}{\gamma_2}\}$ 
for graphs with diameter at least 4. As we have noticed above, the folded Johnson graph with intersection array $\{36,25,16;1,4,18\}$ and the folded halved cube with intersection array $\{66,45,28;1,6,30\}$ do not satisfy (\ref{graph2condition}), however, Theorem \ref{MainTheo2} covers these two cases.

\begin{theo}\label{MainTheo2}
Let $\Gamma$ be a graph of type $2$ and with diameter $D\ge 3$. 
Suppose that $\{x,y,D\}=\{\frac{t-1}{2},\frac{t}{2},t-1+\frac{2}{\gamma_2}\}$.
\begin{itemize}
\item[(i)] If $t\in \{2D,2D+1\}$, then $\gamma_2=2$ and $\Gamma$ is the antipodal quotient of the Johnson 
graph $J(4D,2D)$ or $J(4D+2,2D+1)$, or $\gamma_2=4$ and $\Gamma$ is the antipodal quotient of 
$\frac{1}{2}H(4D,2)$ or $\frac{1}{2}H(4D+2,2)$.
\item[(ii)] If $t=D+1-\frac{2}{\gamma_2}$, then $\gamma_2=4$ 
and $\Gamma$ is the halved graph $\frac{1}{2}H(2D+1,2)$ of the $(2D+1)$-cube.
\end{itemize}
\end{theo}

In \cite{Terw86}, Terwilliger shows that $D\ge 14$ implies that (\ref{graph2condition}) holds, 
which in its turn yields Theorem \ref{Terw86theo}. We will make use of the results from 
the previous sections in order to show that the assumptions on $x,y,t$ and $D$ in Theorem \ref{MainTheo2} 
imply that the local graphs of $\Gamma$ are strongly regular with smallest eigenvalue $-2$, and 
by \cite{BlokBrouwer84}, \cite{BussNeu}, \cite{Neu85} this information is sufficient to determine $\Gamma$. 

\begin{lemma}\label{h2gamma}
Suppose that $\{x,y,D\}=\{\frac{t-1}{2},\frac{t}{2},t-1+\frac{2}{\gamma_2}\}$. 
Then $h=2\gamma_2$.
\end{lemma}
\proof Substituting $x,y,D$ into (\ref{type2def1}) with $i=1$, and taking 
into account that $c_1=1$, we get $h=2\gamma_2$.\wbull

\begin{lemma}$(\cite[Corollary~4.12]{SubAlgPaper})$\label{type2roots}
The Terwilliger polynomial of a graph of type $2$ has the following roots:
\begin{eqnarray}\label{type2roots1}
-1-\frac{(x-1)(D-1)(t-1)}{(x-t+1)(D-t+1)(t-3)}=-1-\frac{b_1(y-t+1)}{(t-1)(y-1)},\\
-1-\frac{(y-1)(D-1)(t-1)}{(y-t+1)(D-t+1)(t-3)}=-1-\frac{b_1(x-t+1)}{(t-1)(x-1)},\\
-1-\frac{(x-1)(y-1)(t-1)}{(x-t+1)(y-t+1)(t-3)}=-1-\frac{b_1(D-t+1)}{(t-1)(D-1)}=\hat{\theta}_D,\label{type2roots11}
\end{eqnarray}
\begin{equation}\label{type2roots2}
-1-\frac{1-t}{3-t}=-2-\frac{2}{t-3}=\hat{\theta}_1.
\end{equation}
\end{lemma}

\begin{lemma}\label{type2leadingterm}
The leading term coefficient of the Terwilliger polynomial of a type $2$ graph with diameter $D\geq 3$ is 
equal to:
\[-p^1_{23}\tau_0=-\frac{2p^1_{23}}{b_1}\times \frac{(t-3)^2}{(t-2)(t-5)}.\]
\end{lemma} 
\proof
It follows directly from Lemma \ref{pplpl}, Lemma \ref{pmimi}, 
and $\theta^*_i=b^*_0+h^*i(i-t)$, $0\leq i\leq D$.\wbull
\medskip

To the rest of this section, we suppose that $\Gamma$ is a type 2 graph 
with diameter $D\geq 3$ and its intersection numbers defined by (\ref{type2def1})--(\ref{type2def2}). 

\begin{lemma}\label{t1D}
If $t=D+1-\frac{2}{\gamma_2}$, $\{x,y\}=\{\frac{t}{2},\frac{t-1}{2}\}$ then $\Gamma$ 
is the halved cube $\frac{1}{2}H(2D+1,2)$ (and $\gamma_2=4$).
\end{lemma}
\proof We first note that if $t=D+1-\frac{2}{\gamma_2}$ and $\{x,y\}=\{\frac{t}{2},\frac{t-1}{2}\}$ 
then $t^*=D+\frac{1}{2}$ by (\ref{ttxyd}), and $h=2\gamma_2$ by Lemma \ref{h2gamma}. 
By (\ref{type2theta}), we have 
$\theta_0>\theta_D>\theta_1>\theta_{D-1}>\theta_2>\ldots$, and hence 
$\Gamma$ has the following $Q$-polynomial ordering: 
\begin{equation}\label{t1DQordering}
E_0,E_2,E_4,\ldots,E_3,E_1,
\end{equation}
where $E_0,E_1,\ldots,E_D$ is the natural ordering of primitive idempotents.

The second largest eigenvalue of $\Gamma$ is $\theta_{D}$, 
defined by (\ref{type2theta}), i.e., 
$\theta_D$ equals $b_0+hD(D-t^*)=b_0-Dh/2=b_0-D\gamma_2$.

Substituting $t=D+1-\frac{2}{\gamma_2}$ and $h=2\gamma_2$ into (\ref{type2defbi}) with $i=1$, 
we see that $\theta_D=b_1-1$. By \cite[Theorem 4.4.11]{BCN}, 
$\Gamma$ is one of the following: 
\begin{itemize}
\item a Hamming graph, a Doob graph, or locally Petersen graph 
(the Doro graph or the Conway -- Smith graph) (and $c_2=2$),
\item a Johnson graph (and $c_2=4$),
\item a halved cube (and $c_2=6$),
\item the Gosset graph with intersection array $\{27,10,1;1,10,27\}$ (and $c_2=10$).
\end{itemize}

The Gosset graph (with respect to the second largest eigenvalue) and the Johnson graphs 
are $Q$-polynomial with type 2A. Since $\gamma_2\ge 1$, we have $c_2>2$ 
and, hence, $\Gamma$ is the halved cube. 
The halved cube $\frac{1}{2}H(m,2)$ has the second 
$Q$-polynomial ordering as in (\ref{t1DQordering}) if and only if $m=2D+1$, see \cite{BI}.
The lemma is proved.\wbull
\medskip

For the remaining cases $t\in \{2D,2D+1\}$, 
we need the following lemma.

\begin{lemma}\label{Interval}
Let $\Delta$ be a non-complete graph on $v$ vertices, 
regular with valency $k$. Let $\eta_1=k,\eta_2,\dots,\eta_v$ be 
all the eigenvalues of $\Delta$. Suppose that there are real numbers $r,s$, $r<s$, 
such that, for every $i=2,\dots,v$, $\eta_i\not\in (r,s)$ holds.  
Then 
\begin{equation}\label{IntervalMain}
kv-k^2+k(r+s)+(v-1)rs\ge 0
\end{equation}
with equality if and only if $\Delta$ is a 
strongly regular graph with non-principal eigenvalues $r,s$.
\end{lemma}
\proof Suppose that, for every $i=2,\dots,v$, $\eta_i\le r$ or 
$\eta_i\ge s$ holds. Since $(\eta_i-r)(\eta_i-s)\ge 0$ holds for every $i>1$, 
it follows that $\sum_{i=2}^v(\eta_i-r)(\eta_i-s)\ge 0$ with equality if and 
only if $\eta_i\in \{r,s\}$ for all $i$, and in this case $\Delta$ is strongly regular. 
Further, we have 
\begin{equation}\label{Interval1}
\sum_{i=2}^v(\eta_i-r)(\eta_i-s)=\sum_{i=2}^{v}\eta_i^2-(r+s)\sum_{i=2}^v\eta_i+(v-1)rs\geq 0.
\end{equation}

Let $B$ be the adjacency matrix of $\Delta$. Then 
\begin{equation}\label{Interval2}
tr(B)=k+\sum_{i=2}^v\eta_i=0,~~tr(B^2)=k^2+\sum_{i=2}^{v}\eta_i^2=kv.
\end{equation}

Combining (\ref{Interval1}) and (\ref{Interval2}), we obtain the required inequality (\ref{IntervalMain}).\wbull

\begin{lemma}\label{t2D}
If $\{x,y,D\}=\{\frac{t-1}{2},\frac{t}{2},t-1+\frac{2}{\gamma_2}\}$ 
and $t\in \{2D,2D+1\}$ then the local graph of any vertex of $\Gamma$ is a strongly regular graph 
with smallest eigenvalue $-2$ and parameters 
\begin{equation}\label{type2srg}
(k,\lambda,\mu)=((t-1)\gamma_2, t\gamma_2/2-2,\gamma_2).
\end{equation}
\end{lemma}
\proof We note that the first three roots (\ref{type2roots1})--(\ref{type2roots11})
of the Terwilliger polynomial are invariant with respect to a permutation of $x,y,D$. 
Substituting $\{x,y,D\}=\{\frac{t-1}{2},\frac{t}{2},t-1+\frac{2}{\gamma_2}\}$ into 
(\ref{type2roots1})--(\ref{type2roots11}), we obtain the following values:
\[
\frac{\gamma_2(t-2)}{2},~~-1+\frac{(\gamma_2(t-2)+2)(t-1)}{2(t-3)},~~\text{and~~}-2,
\]
so that
\[
-2-\frac{2}{t-3}<-2<\frac{\gamma_2(t-2)}{2}<-1+\frac{(\gamma_2(t-2)+2)(t-1)}{2(t-3)}.
\]

If $t\in \{2D,2D+1\}$ then, by Lemma \ref{type2leadingterm}, we see that 
the leading term coefficient of the Terwilliger polynomial is negative, 
and, by Theorem \ref{TerwPolyTheo}, every non-principal eigenvalue $\eta$ of 
the local graph of an arbitrary vertex of $\Gamma$ 
satisfies $\eta\notin (-2,\frac{\gamma_2(t-2)}{2})$.

As $\{x,y,D\}=\{\frac{t-1}{2},\frac{t}{2},t-1+\frac{2}{\gamma_2}\}$ and $t\in \{2D,2D+1\}$, 
we have $h=2\gamma_2$ by Lemma \ref{h2gamma}, and, further, 
\[b_0=\frac{hxyD}{t-1}=t+\gamma_2t(t-1)/2,\]
\[a_1=b_0-b_1-1=\gamma_2(t-1).\]

Applying Lemma \ref{Interval} with $k=a_1$, $v=b_0$, 
$r=-2$, and $s=\frac{\gamma_2(t-2)}{2}$, we get equality in (\ref{IntervalMain}), and hence 
the local graph of any vertex of $\Gamma$ is a strongly regular graph with non-principal eigenvalues 
$-2$ and $\frac{\gamma_2(t-2)}{2}$, 
with parameters (\ref{type2srg}), which shows the lemma.\wbull
\medskip

We recall Seidel's classification of strongly regular graphs with smallest eigenvalue $-2$, 
see \cite{Seidel}.

\begin{theo}\label{Classification2}
A strongly regular graph with smallest eigenvalue $-2$ is one of the following:
\begin{itemize}
\item[(1)] the complete multipartite graph $K_{m\times 2}$ with $m$ parts, each of size $2$, 
with parameters $(2m,2m-2,2m-4,2m-2)$, $m\ge 2$,
\item[(2)] the $m\times m$-grid with parameters $(m^2,2(m-1),m-2,2)$, $m\ge 3$,
\item[(3)] the Shrikhande graph with parameters of $(16,6,2,2)$,
\item[(4)] the triangular graph $T(m)$ with parameters $({m\choose 2},2(m-2),m-2,4)$, $m\ge 5$,
\item[(5)] the three Chang graphs with parameters of $(28,12,6,4)$,
\item[(6)] the Petersen graph with parameters $(10,3,0,1)$,
\item[(7)] the Clebsch graph with parameters $(16,10,6,6)$,
\item[(8)] the Schl\"{a}fli graph with parameters $(27,16,10,8)$.
\end{itemize}
\end{theo}

Finally, we examine which of graphs from Theorem \ref{Classification2} may appear as 
the local graphs of $\Gamma$ if $t\in \{2D,2D+1\}$, and thereby we determine $\Gamma$. 

\begin{lemma}\label{finallemma}
If $\{x,y,D\}=\{\frac{t-1}{2},\frac{t}{2},t-1+\frac{2}{\gamma_2}\}$ 
and $t\in \{2D,2D+1\}$ then $\Gamma$ is the folded Johnson graph 
$\widetilde J(2t,t)$ (and $\gamma_2=2$) or the folded halved cube 
$\frac{1}{2}\widetilde H(2t,2)$ (and $\gamma_2=4$).
\end{lemma}
\proof If $t\in \{2D,2D+1\}$ then, by Lemma \ref{t2D}, 
for every vertex $x\in \Gamma$, $\Gamma(x)$ is a strongly regular graph 
with smallest eigenvalue $-2$ and parameters defined by (\ref{type2srg}).

If $\Gamma(x)$ is the complete multipartite graph with parameters $(2m,2m-2,2m-4,2m-2)$ 
then $k=\mu$, i.e., $(t-1)\gamma_2=\gamma_2$. As clearly $m>2$ and hence 
$\gamma_2>0$ holds, we have $t=2$, a contradiction.

If $\Gamma(x)$ is the $m\times m$-grid (or the Shrikhande graph with parameters 
of $4\times 4$-grid) with parameters $(m^2,2(m-1),m-2,2)$ then $\gamma_2=2$, $t=m$, 
and $c_2=4$. Since $D\geq 3$, we have $t=m\geq 6$, and hence $\Gamma(x)$ is 
the $m\times m$-grid. It is easily seen that the subgraph induced by $\Gamma(u)\cap \Gamma(w)$ 
for a pair of vertices $u,w\in V(\Gamma)$ with $d(u,w)=2$ is a 4-cycle. 
Therefore, by \cite[Theorem~1]{BlokBrouwer84}, $\Gamma$ 
is the antipodal quotient of the Johnson graph $J(2m,m)$.
We note that by \cite{MakPad}, there exist exactly two locally Shrikhande graphs, 
however, both are not distance-regular.

If $\Gamma(x)$ is the triangular graph $T(m)$ (or one of the three Chang graphs with 
parameters of $T(8)$) with parameters $({m\choose 2},2(m-2),m-2,4)$ then $\gamma_2=4$ and $2t=m$.
Since $D\geq 3$, we have $m\geq 12$, and $\Gamma(x)$ is the triangular graph $T(m)$. 
Further, $\Gamma$ has the same intersection array as the antipodal quotient of the halved cube 
$\frac{1}{2}H(2t,2)$. By Step 11 in the proof of \cite[Theorem~3.3]{BussNeu}, 
$\Gamma$ is the antipodal quotient of the halved cube $\frac{1}{2}H(2t,2)$.

If $\Gamma(x)$ is the Petersen graph with parameters $(10,3,0,1)$ 
then $\gamma_2=1$ and $t=4$. This yields that $D=2$, a contradiction.
If $\Gamma(x)$ is the Clebsch graph with parameters $(16,10,6,6)$ 
then $\gamma_2=6$ and $3t-2=6$, $t=8/3$, a contradiction. 
Finally, if $\Gamma(x)$ is the Schl\"{a}fli graph with parameters $(27,16,10,8)$ 
then $\gamma_2=8$ and $4t-2=10$, $t=3$, a contradiction.
The lemma is proved.\wbull
\medskip

Theorem \ref{MainTheo2} now follows from Lemma \ref{t1D} and Lemma \ref{finallemma}. 
Theorem \ref{MainCoro2} follows from Theorem \ref{MainTheo2} and Remark \ref{knowntype2}.

\section*{Acknowledgment}
Part of this work was done while the first author was visiting Tohoku University as a JSPS Postdoctoral Fellow. 
ALG is also supported by the Grant of the President of Russian Federation for young scientists 
(pr. MK-1719.2013.1) and the RFBR grant (pr. 12-01-31098).
JHK thanks for the support of the '100 talents program' of the Chinese Academy of Sciences. 
He also would like to thank Paul Terwilliger for reminding him in June 2011 at the Bled conference 
of the existence of the Terwilliger polynomial.

\end{document}